\documentclass[DIV=calc, paper=letter, fontsize=11pt]{scrartcl}	 

\usepackage[]{units}
\usepackage{lipsum} 
\usepackage[english]{babel} 
\usepackage[protrusion=true,expansion=true]{microtype} 
\usepackage{amsmath,amssymb,amsfonts,amsthm} 
\usepackage[svgnames, table]{xcolor} 
\usepackage[hang, small,labelfont=bf,up,textfont=it,up]{caption} 
\usepackage{booktabs} 
\usepackage{fix-cm}	 
\usepackage[]{units}

\usepackage{pifont} 

\usepackage{tabularx} 

\usepackage{hyperref}
\hypersetup{
colorlinks=true, 
linkcolor=blue, 
citecolor=blue, 
filecolor=blue, 
runcolor=blue, 
urlcolor = blue
}

\usepackage{fullpage}
\usepackage{sectsty} 

\usepackage{fancyhdr} 
\usepackage{lastpage} 

\usepackage{nicefrac}
\usepackage{cite}

\usepackage{tikz}
\usetikzlibrary{calc,patterns,decorations.pathmorphing,decorations.markings}
\usepackage{pgfplots}

\usepackage{multicol}
\usepackage{float} 

\usepackage{subcaption}

\newcounter{tempcolnum}
\makeatletter
\newcommand{\multicolinterrupt}[1]{
\setcounter{tempcolnum}{\col@number}
\end{multicols}
#1%
\begin{multicols}{\value{tempcolnum}}
}
\makeatother

\definecolor{issuePJA_color}{rgb}{1.0,0.0,0.0}

\definecolor{commentPJA_color}{rgb}{1.0,0.0,0.8}

\definecolor{rev_color}{rgb}{0.6,0.0,0.0}

\newcommand{\aHodgeStar}{\overline{\star}}
\newcommand{\aExtD}{\overline{\mb{d}}}

\usepackage{graphicx}
\DeclareGraphicsExtensions{.png}

\usepackage{lettrine} 
\newcommand{\initial}[1]{ 
\lettrine[lines=3,lhang=0.3,nindent=0em]{
\color{DarkGoldenrod}
{\textsf{#1}}}{}}

\usepackage{caption}

\usepackage{graphicx}

\usepackage{titling} 

\newcommand{\HorRule}{\color{DarkGoldenrod} \rule{\linewidth}{1pt}} 

\pretitle{\vspace{-80pt} \begin{flushleft} \HorRule \fontsize{20}{20} \color{DarkRed} \selectfont} 

\title{Spectral Numerical Exterior Calculus Methods for Differential Equations on Radial Manifolds} 

\posttitle{\par\end{flushleft}} 
\preauthor{\vspace{-5pt} \begin{flushleft}\fontsize{14}{14} \large  \color{DarkRed}} 
\author{B. Gross and P. J. Atzberger$^{*}$ } 
\postauthor{\fontsize{12}{12} \small \color{Black} 
Department of Mathematics, Department of Mechanical Engineering, University of California Santa Barbara; atzberg@gmail.com; http://atzberger.org/ \\
\vskip 1.5em
\lfoot{\footnotesize * Work supported by DOE ASCR CM4 DE-SC0009254, NSF CAREER Grant DMS-0956210, and \\ NSF Grant DMS - 1616353.}
\end{flushleft} \vspace{-18pt} \HorRule} 

\newcommand{\mb}[1]{\mathbf{#1}}
\newcommand{\bs}[1]{\boldsymbol{#1}}

\newcommand{\mDiv}{ {\mbox{div}} }
\newcommand{\mGrad}{ {\mbox{grad}} }
\newcommand{\mCurl}{ {\mbox{curl}} }

\DeclareGraphicsExtensions{.png}

\begin{document}

\maketitle 

\thispagestyle{fancy} 


\initial{W}\textbf{e develop exterior calculus approaches for partial differential equations on radial manifolds.   We introduce numerical methods that approximate with spectral accuracy the exterior derivative $\mb{d}$, Hodge star $\star$, and their compositions. To achieve discretizations with high precision and symmetry, we develop hyperinterpolation methods based on spherical harmonics and Lebedev quadrature.  We perform convergence studies of our numerical exterior derivative operator $\aExtD$ and Hodge star operator $\aHodgeStar$ showing each converge spectrally to $\mb{d}$ and $\star$.  We show how the numerical operators can be naturally composed to formulate general numerical approximations for solving differential equations on manifolds.  We present results for the Laplace-Beltrami equations demonstrating our approach.}

\setlength{\parindent}{5ex}

\section{Introduction}
\label{sec:Intro}
There has been a lot of recent interest in numerical methods related to exterior calculus~\cite{Hirani2003,DesbrunSiggraph2013, Arnold2010,HolstExCalcFEM2017}.  Application areas include hydrodynamics within fluid interfaces~\cite{AtzbergerSoftMatter2016, ArroyoRelaxationDynamics2009}, electrodynamics~\cite{MarsdenElectrodynamics2015}, stable methods for finite elements~\cite{Arnold2006,Arnold2010,HolstExCalcFEM2017}, and geometric processing in computer graphics~\cite{Desbrun2003,DesbrunSiggraph2013,DesbrunSubdivision2016,Schroder2006}.  The exterior calculus of differential geometry provides a cooordinate invariant way to formulate equations on manifolds with close connections to topological and geometric structures inherent in mechanics~\cite{Marsden1994,Marsden2007}.  The exterior calculus provides less coordinate-centric expressions for analysis and numerical approximation that often can be interpreted more readily in terms of the geometry than alternative approaches such as the tensor calculus~\cite{Eells1955,Whitney1957,Marsden2007}.  Many discrete exterior calculus approaches have been developed for efficient low-order approximations on triangulated meshes~\cite{Hirani2003, Desbrun2003, DesbrunSiggraph2013}.  There has also been recent work on higher-order methods based on collocation for product bases~\cite{DesbrunSpectralChainCollocation2014}, finite element differential forms~\cite{Arnold2006}, and more recently subdivision surfaces~\cite{DesbrunSubdivision2016}.  In these discrete exterior calculus approaches an effort is made to introduce approximations for fundamental operators such as the exterior derivative $\mb{d}$ and Hodge star $\star$ operators that on the discrete level preserve inherent geometric relations~\cite{DesbrunSiggraph2013,Arnold2006,HolstExCalcFEM2017}.  These operators are then used through composition to perform geometric processing tasks or approximate differential equations.  We show here how related approaches can be developed using hyperinterpolation to obtain spectrally accurate methods for exterior calculus on radial manifolds.

\begin{figure}[H]
\centering
\includegraphics[width=0.7\linewidth]{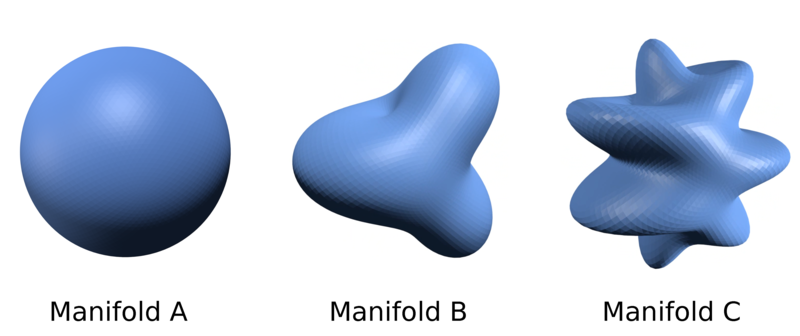}
\caption{Radial Manifolds: A radial manifold is defined as a surface where each point can be connected by a line segment to the origin without intersecting the surface.  Shown are three radial manifolds which for discussions we refer to interchangeably as the (i) Sphere / Manifold A, (ii) Dimple / Manifold B, and (iii) Fountain / Manifold C.  The manifolds are generated by the radial functions (i) $r(\theta,\phi) = 1.0$, (ii) $r(\theta,\phi) = 1 + r_0\sin(3\phi)\cos(\theta)$ with $r_0 = 0.4$, and (iii) $r(\theta,\phi) = 1 + r_0 \sin(7\phi)\cos(\theta)$ with $r_0 = 0.4$.  For additional discussion of the differential geometry of radial manifolds see Appendix~\ref{appendix:radialManifold}.}
\label{fig:RadialManifolds}
\end{figure}

\lfoot{} 

We take an isogeometric approach based on the hyperinterpolation of spherical harmonics to represent the manifold geometry, differential forms, and related scalar and vector fields~\cite{Hughes2009, Sloan1995, HanBookSphericalHarmonics2010}.  Hyperinterpolation methods use an oversampling of functions to grapple with some of the inherent challenges in designing optimal nodes for interpolation on general domains~\cite{Sloan1995,Sloan2001,Sloan2000,Reimer2000}.  This allows for the treatment of approximation instead using approaches such as $L^2$-orthogonal projection based on exact quadratures up to a desired order~\cite{Sloan1995,Sloan2012}.  To achieve discretizations on spherical toplogies with favorable symmetry, we use the nodes of Lebedev quadrature~\cite{Lebedev1999,Lebedev1976}.  While the more common approach of using samping points based on lattitude and longitude does provide some computational advantage through fast transforms, the sampling points have poor symmetry and exhibit a significant inhomogeneous distribution over the surface with many sample points clustering near the poles~\cite{Schaeffer2013,Healy2003,HealyDriscoll1994}.  The Lebedev quadrature points are more regularly distributed on the surface and for a comperable number of points provide over the surface a more uniform resolution of functions.  Furthermore, the Lebedev quadrature points are invariant under rotations corresponding to octohedral symmetry~\cite{Lebedev1999,Lebedev1976}.  We use an $L^2$-projection to spherical harmonics to approximate the exterior derivative operator $\mb{d}$ and Hodge star operator $\star$ on the surface.  We show that our methods provide spectrally accurate approximations for these operators and their compositions.  We show how the methods can be used to develop numerical solvers for differential equations on manifolds.  We present results for the Laplace-Beltrami equations demonstrating our approach.  The introduced numerical methods can be applied quite generally for approximating with spectral accuracy differential equations or other exterior calculus operations on radial manifolds.

\section{Differential geometry and conservation laws on manifolds}

We briefly discuss the formulation of conservation laws on manifolds in covariant form and discuss related concepts in differential geometry.  More detailed discussions of the associated differential geometry can be found in~\cite{Abraham1988,Spivak1971,Pressley2001}. Differential forms arise as a natural approach in formulating conservation laws and relations in continuum mechanics~\cite{Eells1955,Marsden2007}.  The exterior calculus provides a convenient means to generalize the Stokes Theorem, Divergence Theorem, and Green's Identities to curved spaces~\cite{Marsden1994}.  We introduce numerical methods for approximating exterior calculus operations on radial manifolds.  A radial manifold is defined as a surface where every point can be connected by a line segment to the origin without intersecting the surface.  We consider compact radial manifolds without boundaries which correspond to shapes with spherical topology, see Figure~\ref{fig:RadialManifolds}.

\begin{figure}[H]
\centering
\includegraphics[width=0.8\linewidth]{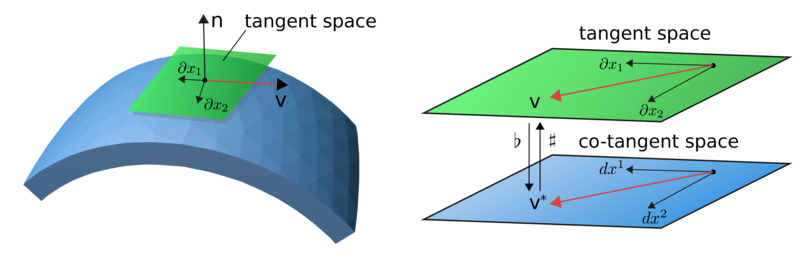}
\caption{We consider vectors and covectors on the manifold.  The isomorphisms between the tangent space and co-tangent space are given by ${\flat}[\mb{v}] = v^j \partial_{\mb{x}^j} \rightarrow v_i d\mb{x}^i = \mb{v}^*$ and ${\sharp}[ \mb{v}^*] = v_i d\mb{x}^i \rightarrow v^j\partial_{\mb{x}^j} = \mb{v}$.  The musical notation $\sharp$ and $\flat$ is used to indicate how the indices are either being raised or lowered in the tensor notation~\cite{Abraham1988,Spivak1971}.}
\label{fig:cotangentSpace}
\end{figure}

\subsection{Covariant Form and Exterior Calculus}
We briefly discuss the covariant formulation of tensors and exterior calculus.  More detailed discussions can be found in~\cite{Abraham1988,Spivak1971}.  The covariant formulation relies on casting relationships using the cotangent bundle of the manifold.  The cotangent space $\mathcal{V}^*$ at location $\mb{x}$ is a vector space that consists of all of the linear functionals that act on the tangent space $\mathcal{V}$.  The dual vector $\mb{v}^{*} \in \mathcal{V}^*$ acts on a vector $\mb{u} \in \mathcal{V}$ as $\mb{v}^{*}[\mb{u}] = \langle \mb{v},\mb{u} \rangle = v^i g_{ij} u^j = v_j u^j$ where $g_{ij}$ denotes the entries of the metric tensor of the manifold~\cite{Pressley2001,Abraham1988}.  We use the lower index notational conventions for $v_j =  g_{ij} v^{i}$ to obtain from the tangent vector components $v^{i}$ the associated covariant vector components $v_j$. The Einstein summation conventions are used throughout.  We adopt the usual convention of indexing contravariant tensors with superscripts and covariant tensors with subscripts~\cite{Abraham1988,Spivak1971}.  Since each dual element $\mb{v}^*$ can be expressed through an inner-product with a representative tangent vector $\mb{v}$, the dual space $\mathcal{V}^*$ is isomorphic to the tangent space $\mathcal{V}$.  We adopt the musical notation $\sharp$ and $\flat$ for raising and lowering indices.  This corresponds to the isomorphisms between the tangent and co-tangent spaces of the surface given by
\begin{eqnarray}
\label{equ_iso}
{\flat} :   v^j \partial_{\mb{x}^j} \rightarrow 
v_i d\mb{x}^i, \hspace{1cm}
{\sharp} :  v_i d\mb{x}^i \rightarrow v^j \partial_{\mb{x}^j}.
\end{eqnarray}
The $\partial_{\mb{x}^j}$ denotes for a given choice of coordinates ${x}^i$ the associated basis vectors for the tangent space $\mathcal{V}$.  The $d{\mb{x}^j}$ denotes for a given choice of coordinates the $1$-form basis elements for the co-tangent space $\mathcal{V}^*$~\cite{Abraham1988}.  The isomorphisms can also be expressed directly in terms of the components as $v_i = g_{ij} v^j$ and $v^i = g^{ij} v_j$, where we denote the metric tensor as $g_{ij}$ and its inverse as $g^{ij}$.  These conventions extend naturally to higher rank tensors~\cite{Abraham1988}.

We consider the exterior derivative $\mb{d}$ and the Hodge star $\star$ on the manifold.  The exterior derivative $\mb{d}$ acts on a $k$-form $\alpha$ to yield the $k+1$-form 
\begin{eqnarray}
\mb{d} \alpha = 
\frac{1}{k!}
\frac{\partial \alpha_{i_1 \ldots i_k}}{\partial x^j} \mb{d}\mb{x}^j \wedge  \mb{d}\mb{x}^{i_1} \wedge \cdots \wedge \mb{d}\mb{x}^{i_k}.
\end{eqnarray}
Here we have taken that the $k$-form to be $\alpha = \left(1/k!\right)\alpha_{i_1 \ldots i_k}\mb{d}\mb{x}^{i_1} \wedge \cdots \wedge \mb{d}\mb{x}^{i_{k}}$, where $\wedge$ denotes the wedge product~\cite{Spivak1971,Abraham1988}.  Since we are using the Einstein summation conventions, all permutations of the indices can arise and hence the factor of $1/k!$.  For notational convenience and to make expressions more compact we will sometimes adsorb this factor implicity into $\alpha_{i_1 \ldots i_k}$.  The Hodge star $\star$ acts on a $k$-form and yields the $(n-k)$-form
\begin{eqnarray}
\label{eqn_HodgeStar}
\star \alpha = 
\frac{\sqrt{|g|}}{(n-k)!k!}
\alpha^{i_1 \ldots i_k}\epsilon_{i_1\ldots i_k j_1 \ldots j_{n-k}} 
\cdot \mb{d}\mb{x}^{j_1} \wedge \cdots \wedge \mb{d}\mb{x}^{j_{n-k}},
\end{eqnarray}
where 
$\alpha^{i_1 \ldots i_k} = g^{i_1 \ell_1}\cdots g^{i_k \ell_k} \alpha_{\ell_1 \ldots \ell_k}$,
$\sqrt{|g|}$ is the square-root of the determinant of the metric tensor, and 
$\epsilon_{i_1\ldots i_k j_1 \ldots j_{n-k}}$ denotes the 
Levi-Civita tensor~\cite{Marsden1994}.  The Hodge star $\star$ can be thought of intuitively as yielding the orthogonal compliment to the $k$-form in the sense that the wedge product recovers a scaled volume form~\cite{Hirani2003,DesbrunSiggraph2013}.  This interpretation of the Hodge star is easiet to see when the local coordinates are chosen to be orthonormal at location $\mb{x}$~\cite{Spivak1971,Abraham1988}.

Many of the operations in vector calculus can be generalized to manfolds using the exterior calculus as
\begin{eqnarray}
\label{equ:gradDivCur}
\mGrad(f) = \lbrack \mb{d}f\rbrack^{\sharp}, \hspace{1cm}
\mDiv(\mb{F}) =  -(\star \mb{d}\star \mb{F}^\flat) = -\bs{\delta} \mb{F}^\flat, \hspace{1cm}
\mCurl(\mb{F}) = \left\lbrack\star (\mb{d}\mb{F}^\flat) \right\rbrack^{\sharp},
\end{eqnarray}
where $f$ is a smooth scalar function and the $\mb{F}$ is a smooth vector field.  We use the notation $\bs{\delta}$ to denote the co-differential operator given by $\bs{\delta} = \star \mb{d} \star$.  As can be seen above, the $\bs{\delta}$ operator is closely related to the divergence operation in vector calculus.  It is also natural to consider common vector calculus differential operators such as the Laplacian.  It is important to note that on manifolds there are a few different operators that share many of the features with the Laplacian of vector calculus.  This requires care when formulating conservation laws or considering constitutive models.  A few generalizations of the Laplacian include 
\begin{eqnarray}
\Delta^H (\mb{F}) 
 =  -\left\lbrack \left(\bs{\delta}\mb{d} + \mb{d}\bs{\delta} \right) \mb{F}^\flat\right\rbrack^{\sharp},\hspace{0.4cm}
\Delta^S (\mb{F}) 
 =  -\left\lbrack \bs{\delta}\mb{d}  \mb{F}^\flat\right\rbrack^{\sharp},\hspace{0.4cm}
\Delta^H f  =  \Delta^R f = -(\star \mb{d}\star )\mb{d} f = -\bs{\delta} \mb{d} f. 
\end{eqnarray}
The $\Delta^R = \mDiv(\mGrad(\cdot))$ denotes the rough-Laplacian given by the usual divergence of the gradient.  For vector fields, $\Delta^H (\mb{F})$ denotes the Hodge-de Rham Laplacian, which has similarities to taking the curl of the curl~\cite{Abraham1988}.  In fact, in the case that $\mDiv(\mb{F}) = -\bs{\delta} \mb{F}^{\flat} = 0$ we have $\Delta^H (\mb{F}) = \Delta^S (\mb{F}) = -\left[\bs{\delta} \mb{d} \mb{F}^{\flat}\right]^{\sharp}$.

The exterior calculus provides generalizations of the Stokes Theorem and the Divergence Theorem to manifolds.  These take the form respectively
\begin{eqnarray}
\label{equ:exteriorStokesDivThms}
\int_{\partial \Omega} \omega = \int_{\Omega} \mb{d}\omega, \hspace{1cm} \int_{\partial \Omega} \star \omega = \int_{\Omega} \mb{d} \star \omega.
\end{eqnarray}
The $\omega$ is a $k$-form and $\Omega$ denotes a general smooth domain with smooth boundary $\partial \Omega$ within the manifold.  When $\Omega$ is in two dimensions and $\omega$ is a $1$-form, the integrals on the left-hand-side perform a type of line integral over the boundary contour $\partial \Omega$.  For the Stokes Theorem the component of the vector field that is tangent to the contour is integrated.  For the Divergence Theorem the Hodge star $\star$ ensures that only the component of the vector field that is normal to the contour is integrated.  The right-hand side then integrates over the interior region of $\Omega$.  For the Stokes Theorem $\mb{d} \omega$ corresponds to a generalized representation of the curl operation.  For the Divergence Theorem $\mb{d}\star \omega$ provides a representation of the divergence.  Additional applications of the Hodge star $\star$ and isomorphisms $\flat$, $\sharp$ can also be useful to make conversions that bring these expressions into closer agreement with the standard vector calculus interpretations and intuition.  Notice the relation of these expressions when $\omega = \mb{F}^{\flat}$ to equation~\ref{equ:gradDivCur}.  A useful feature of equation~\ref{equ:exteriorStokesDivThms} is that it also generalizes readily to higher dimensions and $k$-forms.

\subsection{Conservation Laws on Manifolds}
For conservation laws, the exterior caclulus provides convenient ways to formulate relations and constitutive laws without relying on cumbersome coordinate expressions.  One application of this idea can be seen by building on the generalized Stokes Theorem and the generalized Divergence Theorem in equation~\ref{equ:exteriorStokesDivThms}.  For a conserved scalar quantity $u$, let $\omega$ denote a differential $1$-form that represents the local flux over a boundary, such as mass or energy.  Conservation of $u$ and application of the generalized Divergence Theorem gives 
\begin{eqnarray}
\label{equ:consrvLaw}
\frac{\partial}{\partial t} \int_{\Omega} \star u = \int_{\partial \Omega} \star \omega = \int_{\Omega} \mb{d} \star \omega.
\end{eqnarray}
Since $u$ is a scalar field on the manifold ($0$-form), we represent its integral over the surface as the $2$-form $\star u = \tilde{u}_{ij} d\mb{x}^i\wedge d\mb{x}^j$.  The flux is represented by the $1$-form $\omega = \omega_i d\mb{x}^i$.  By the generalized Divergence Theorem, we have that the exterior derivative $\mb{d}$ gives the $2$-form $\mb{d} \star \omega = \tilde{\omega}_{ij} d\mb{x}^i\wedge d\mb{x}^j$ that integrates over $\Omega$ to give the same value as the flux $\star \omega$ integrated over the boundary $\partial \Omega$.  Since $\Omega$ is arbitrary, this provides for smooth fields the local expression for the conservation law in covariant form
\begin{eqnarray}
\frac{\partial u}{\partial t} = -\star \mb{d} \star \omega = -\bs{\delta} \omega.
\end{eqnarray}
To obtain this result, we applied the Hodge star $\star$ to both sides of equation~\ref{equ:consrvLaw} and used that 
$\star\star = -1^m$, where $m = k\cdot(n-k)$ for a $k$-form in $n$ dimensional space~\cite{Abraham1988,Spivak1971}.  The co-differential is given by $\bs{\delta} = \star \mb{d} \star$ and we see it plays here a role very similar to the vector calculus operation of divergence.

In many cases the flux itself follows a law depending on the conserved quantity $u$.  In the case of Fourier's Law for heat conduction or Fick's Law for mass diffusion, the flux depends on the gradient of $u$.  This generalizes on the manifold to
\begin{eqnarray}
\omega = \mb{d} u.
\end{eqnarray}
The local conservation law in covariant form is
\begin{eqnarray}
\frac{\partial u}{\partial t} = -\bs{\delta} \mb{d} u.
\end{eqnarray}
When $u$ is heat or mass, this gives the generalization of the heat equation or the diffusion equation to the manifold.  

We have given a brief discussion of the utility of exterior calculus methods in formulating scalar conservation laws on manifolds.  Similar approaches can also be taken for vector conservation laws on manifolds and in the formulation of further constituitive laws and relations in continuum mechanics~\cite{AtzbergerSoftMatter2016, Marsden2007}.  A particular advantage of the exterior calculus is that we can formulate models and constuitive relations following closely intuition developed in the context of vector calculus without the need to deal with cumbersome coordinate expressions.  We see that formulations rely on a composition of the exterior derivative $\mb{d}$ and Hodge star $\star$ operations.  To utilize this approach in numerical calculations we need accurate approximations of the action of these operators on differential forms.

\section{Numerical Methods for Exterior Calculus}
\label{sec:num_methods}
We take an isogeometric approach based on the hyperinterpolation of spherical harmonics to represent the manifold geometry, differential forms, and related scalar and vector fields~\cite{Hughes2009, Sloan1995, HanBookSphericalHarmonics2010}.  In hyperinterpolation, functions are oversampled to avoid many of the inherent issues associated with trying to design an optimal collection of nodes for Lagrange interpolation~\cite{Sloan1995,Sloan2001,Sloan2000,Reimer2000}.  This allows for functions to be approximated through $L^2$-orthogonal projections using exact quadratures up to a desired order~\cite{Sloan1995,Sloan2012}.  To achieve discretizations with favorable symmetry on the sphere, we use the nodes of Lebedev quadrature~\cite{Lebedev1999,Lebedev1976}.  This is in contrast to the more common approach of using samping points based on lattitude and longitude coordinates.  While lattitude-longitude samplings have a computational advantage through fast transforms, the sampling points have poor symmetry and inhomogeneous distribution over the surface with many  points clustering near the poles~\cite{Schaeffer2013,Healy2003,HealyDriscoll1994}.  The Lebedev quadrature points provide a more regular distribution on the surface. For a comperable number of points, the Lebedev sampling provides a more uniform resolution of functions.  The Lebedev quadrature points also have the feature of being invariant under rotations corresponding to octohedral symmetry~\cite{Lebedev1999,Lebedev1976}.  We show Lebedev nodes on example radial manifolds in Figure~\ref{fig:LebedevQuad}.  For additional discussion of the differential geometry of radial manifolds see Appendix~\ref{appendix:radialManifold}.  
\begin{figure}[H]
\centering
\includegraphics[width=0.8\linewidth]{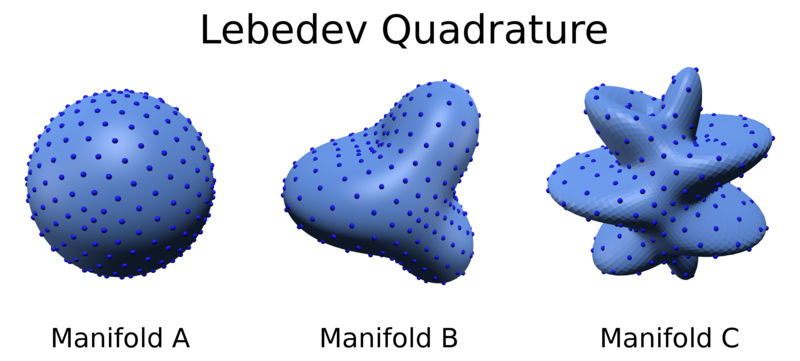}
\caption{Lebedev Quadrature.  Shown are the sample points of the Lebedev quadrature in the case of 302 nodes.  The Lebedev nodes distribute nearly uniformly over the surface and are invariant under the rotations corresponding to octahedral symmetry~\cite{Lebedev1976,Lebedev1999}.}
\label{fig:LebedevQuad}
\end{figure}

\subsection{Hyperinterpolation and $L^2$-Projection}
We use hyperinterpolation to obtain a continuum approximation to fields on the manifold surface~\cite{Sloan1995,Sloan2012}.  To obtain a continuum representation of a function $f$ on the surface, we perform an $L^2$-orthogonal projection $\mathcal{P}$ to the space spanned by spherical harmonics up to order $\lfloor L/2 \rfloor$,
\begin{eqnarray}
\label{equ_proj_op}
\mathcal{P}\left[f \right] = \bar{f}(\theta,\phi) = \sum_i \hat{f}_i Y_i(\theta,\phi),
\end{eqnarray}
where $\hat{f}_i = \langle f, Y_i \rangle_Q$.  We take the spherical harmonics $Y_i$ to be normalized with $\langle Y_j, Y_i \rangle_Q = \delta_{ij}$ see Appendix~\ref{appendix:sphericalHarmonics}.  We use the discrete inner-product defined by
\begin{eqnarray}
\label{equ_innerProd_Q}
\langle u, v \rangle_Q = \sum_\ell w_\ell u(\mb{x}_\ell) v(\mb{x}_\ell).
\end{eqnarray}
where $w_\ell, \mb{x}_\ell$ are the Lebedev quadrature weights and nodes.  When the quadrature is of order $L$ and the functions $u,v$ are each band-limited with respect to spherical harmonics up to order $\lfloor L/2 \rfloor$, the the inner-product is the same as the $L^2$-inner-product $\langle u, v \rangle_Q = \langle u, v \rangle_{L^2}$.  This yields the projection property $\mathcal{P}^2 = \mathcal{P}$, see equation~\ref{equ_proj_op}.  For additional discussion of how we employ the spherical harmonics see Appendix~\ref{appendix:sphericalHarmonics}.  

In practice, computing the inner-product $\langle \cdot, \cdot \rangle_Q$ only requires we know values of $f$ at the Lebedev nodes $\{\mb{x}_\ell\}$.  This is utilized to represent functions on the surface in numerical calculations.  We use this property to represent differential forms on the surface by an equivalent vector field at the Lebedev nodes.  We perform in calculations conversions as needed by using the isomorphisms $\flat$, $\sharp$ in equations~\ref{equ_iso}, see Figure~\ref{fig:cotangentSpace}.  In two dimensions the $0$-forms and $2$-forms on the surface are equivalent to scalar fields (a vector field with only one component).  For the more interesting case of $1$-forms on the surface, we use as our numerical represention an equivalent vector field with values specified at each of the Lebedev quadrature nodes.  To simplify our discussion of our numerical methods we use the terminiology vector field and scalar field interchangably throughout.  

The $1$-form $\mb{v}^{\flat}$ is equivalent through the isomorphism $\sharp$ to the vector field $\mb{v}^{\sharp}$.  We represent $1$-forms by the values of $\mb{v}^{\sharp}$ stored at the Lebedev quadrature nodes $\{\mb{x}_\ell\}$.  In numerical calculations we avoid issues with charts and coordinate singularities by representing the form as an expansion of spherical harmonics using the coordinates of the embedding space.  This is done by representing the components of the associated vector field $\mb{v}^{\sharp}$ using the embedding space basis $\bs{\iota}_1, \bs{\iota}_2, \bs{\iota}_3$ as $\mb{v}^{\sharp}(\mb{x}_\ell) = \bar{v}^{x}\bs{\iota}_1 + \bar{v}^{y}\bs{\iota}_2 + \bar{v}^{z}\bs{\iota}_2 = \left[\bar{v}^{x}, \bar{v}^{y}, \bar{v}^{z}\right]_{\bs{\iota}_1, \bs{\iota}_2, \bs{\iota}_3}$.  Storing the values $\left[\bar{v}^{x}, \bar{v}^{y}, \bar{v}^{z}\right]$ at the Lebedev nodes provides a convenient numerical representation of the differential form.  To simplify the notation, we will often drop the subscript on $[\cdot,\cdot,\cdot]$ for the basis when it can be understood by context.    When a continuum representation of the vector field is needed in our numerical calculations, we use the hyperinterpolation in equation~\ref{equ_proj_op} to obtain the associated smooth vector field
\begin{eqnarray}
\label{equ_v_cont}
\mb{\bar{v}}^{\sharp}(\theta,\phi) = [\mathcal{P}\bar{v}^{x}, \mathcal{P}\bar{v}^{y}, \mathcal{P}\bar{v}^{z}].
\end{eqnarray} 
We take a similar approach for $0$-forms and $2$-forms which are much easier to handle and are represented by the scalar field $\bar{v}$ to yield $\mb{\bar{v}}^{\sharp}(\theta,\phi) = \mathcal{P}\bar{v}$.

\subsection{Numerical Exterior Calculus Operators $\aExtD$ and $\aHodgeStar$}

\label{sec_num_ext}

To approximate the exterior derivative $\mb{d}$, we need to approximate derivatives of our numerical representation at the Lebdev nodes for  differential forms.  For this purpose, we make use of the hyperinterpolation provided by equation~\ref{equ_v_cont}. We remark that our approach in our numerical representation making use of the embedding space basis provides a global description of the differential form over the entire surface and a consistent way to obtain derivatives between different coordinate charts.  For a given chart, a differential form has coordinate components for a $0$-form given by $\mb{v}^{\flat} = v$, a $1$-form by $\mb{v}^{\flat} = v_i d\mb{x}^i$, and a $2$-form by $\mb{v}^{\flat} = v_{ij} d\mb{x}^i\wedge d\mb{x}^j$.  To numerically compute derivatives based on these expressions we perform a conversion from the vector field representation in the embedding space to the local coordinate representation of the components.

The $1$-form presents the most interesting case with the $0$-form and $2$-form handled similarly.  For $1$-forms the components $v_i$ are related to the components $\bar{v}^k$ of the vector field representation by $v_i = g_{ij} v^j = g_{ij} a^j_k \bar{v}^k$.  The term $a^j_k$ converts between the components $\bar{v}^k$ given in the coordinates of the embedding space $\bs{\iota}_1, \bs{\iota}_2, \bs{\iota}_3$ to the components ${v}^j$ given in the local coordinates on the surface $\partial_\theta, \partial_\phi$.  The exterior derivative of a $1$-form can be expressed in coordinate components as $\mb{d}\mb{v}^{\flat} = \partial_s v_i \hspace{0.1cm} d\mb{x}^s \wedge d\mb{x}^i$.  For numerical calculations at a given location $\mb{x}_\ell$ (Lebedev node), we choose an appropriate coordinate chart that is locally non-degenerate.  We compute the component derivative as
\begin{eqnarray}
\label{equ_der_comp_v}
{\partial_s} v_i =  ({\partial_s}g_{ij}) a^j_k \bar{v}^k + g_{ij} ({\partial_s}a^j_k) \bar{v}^k + g_{ij} a^j_k ({\partial_s}\bar{v}^k).
\end{eqnarray}
The first two terms only depend on the geometry of the manifold and only the values of the differential form at location $\mb{x}_\ell$ (Lebedev node).  This can be obtained readily from the spherical harmonics representation of the geometry of the manifold see Appendix~\ref{appendix:sphericalHarmonics} and Appendix~\ref{appendix:radialManifold}.  In contrast, the last term depends on the derivatives of the coordinate components and requires use of the continuum representation from the hyperinterpolation obtained in equation~\ref{equ_v_cont}.  Putting this together with equation~\ref{equ_der_comp_v} and the coordinate expression for the exterior derivative, we obtain a numerical exterior derivative operator $\aExtD$ for $1$-forms.  The case of $0$-form and $2$-form can be handled similarly.  We should mention that the case of a $2$-form in two dimensions has exterior derivative zero which we also impose in our numerical calculations.  In this manner, we obtain a numerical operator $\aExtD$ that maps a $k$-form defined at the Lebedev nodes to a $(k+1$)-form defined at the Lebedev nodes.  This results in a convenient map between our representations useful in compositions for further application of other numerical exterior calculus operations.

We remark that in practical implimentations of the numerical exterior derivative operator $\aExtD$ it is convenient to represent the coordinate conversion in matrix-vector notation as $\mb{v} = G A^{-1} \mb{\bar{v}}$.  The matrix entries $[A^{-1}]_{jk} = a^j_k$ correspond to the change from the coordinates of the embedding space to the local coordinates of the tangent space.  With this notation, the derivatives in local coordinates can be expressed as ${\partial_s} \mb{v} = ({\partial_s}G) A^{-1} \mb{\bar{v}} + G ({\partial_s}A^{-1}) \mb{\bar{v}} + G A^{-1} ({\partial_s}\mb{\bar{v}})$.  To avoid differentiating components of the inverse matrix $A^{-1}$, we use the identity ${\partial_s}A^{-1} = - A^{-1} ({\partial_s}A) A^{-1}$ and use a linear algebra solver to compute the action of $A^{-1}$.  This is done at each Lebedev node $\mb{x}_\ell$ with an appropriate choice made for the  coordinate chart that is locally non-degenerate.  For additional information on how we define the coordinate charts and how we perform practical computations from our representations of the manifold geometry see Appendix~\ref{appendix:radialManifold}.

We approximate next the Hodge star $\star$ operator on differential forms.
We remark that in the related area of discrete exterior calculus (DEC) efforts are made to preserve geometric structure in the discrete setting, often on triangulated meshes.  Interesting issues arise in DEC from the discrete geometry with which one must grapple and extensive studies have been conducted to formulate good approximations for the Hodge star $\star$ operator~\cite{Hirani2016,Hirani2013}.  Here we avoid many of these issues since we treat the operator at the continuum level and have more geometric information available to us from our spectral representation of both the manifold and the differential forms.  

We approximate the Hodge star $\star$ operator on differential forms by a numerical operator $\aHodgeStar$ which makes use of the representation at the Lebedev nodes.  The Hodge star $\star$ has the feature that it is a local operation that involves values of the differential form and metric tensor only at an individual Lebedev node $\mb{x}_\ell$.  We obtain a numerical operator $\aHodgeStar$ by applying the isomorphisms and metric tensor using equation~\ref{equ_iso} and equation~\ref{eqn_HodgeStar}.  The main consideration numerically is to choose well the coordinate chart so it is locally non-degenerate.  The approximation enters through the fidelity of the metric tensor computed from our representation of the manifold geometry.  The geometry for the radial manifold is determined by the radial function $r(\theta,\phi)$ which is represented as an expansion in a finite number of spherical harmonics up to order $L$, $r(\theta,\phi) = \sum_i \hat{r}_i Y_i(\theta,\phi)$.  This provides in a given coordinate chart at location $\mb{x}_\ell$ the metric tensor and curvature tensor along with the local coordinate basis vectors $\partial_\theta, \partial_\phi$ and their derivatives see Appendix~\ref{appendix:radialManifold}.  The numerical Hodge star operator $\aHodgeStar$ maps $k$-forms defined at the Lebedev nodes to $(n-k)$-forms defined at the Lebedev nodes.  This provides a convenient map between representations for further application of numerical exterior calculus operations.  In this manner, the numerical exterior derivative $\aExtD$ operator and numerical Hodge star $\aHodgeStar$ operator can be used through composition to numerically approximate more complex exterior calculus operations on the manifold.

\subsection{Solving Differential Equations on Manifolds}
\label{sec:methodsForPDE}
We consider differential equations on the manifold of the form
\begin{eqnarray}
\label{equ_main_PDE}
\mathcal{L} u = - g.
\end{eqnarray}
The $\mathcal{L}$ denotes a linear differential operator that can be expressed in terms of a composition of the exterior calculus operations $\mb{d}$ and $\star$.  For example, in the case of the Laplace-Beltrami equation this would correspond to the operator $\mathcal{L} = -\bs{\delta} \mb{d} = -\star \mb{d} \star \mb{d}$.  We discretize the operator $\mathcal{L}$ by using a composition of the numerical exterior calculus operators $\aHodgeStar$ and $\aExtD$ to obtain $\mathcal{\tilde{L}}$.  For the Laplace-Beltrami equation, this corresponds to $\mathcal{\tilde{L}} = -\bar{\bs{\delta}} \mb{\bar{d}} = -\aHodgeStar \hspace{0.1cm} \aExtD \hspace{0.1cm} \aHodgeStar \hspace{0.1cm} \aExtD$.  We approximate equation~\ref{equ_main_PDE} on the manifold by
\begin{eqnarray}
\label{equ_main_PDE_weakForm}
\langle \mathcal{\tilde{L}} \bar{u}, Y_i \rangle_Q = - \langle \bar{g}, Y_i \rangle_Q.
\end{eqnarray}
The $\bar{u} = \sum_j \hat{u}_j Y_j$, $\bar{g} = \sum_j \hat{g}_j Y_j$ denote expansions up to order $\lfloor L/2 \rfloor$ and $\langle \cdot,\cdot \rangle_Q$ denotes the Lebedev inner-product computed up to order $L$ in equation~\ref{equ_innerProd_Q}.  This provides a Galerkin approximation where some additional sources of approximation arise from the treatment of the differential operator by $\mathcal{\tilde{L}}$.  The approximation can be expressed in terms of the solution coefficients $\mb{\hat{u}}$ as
\begin{eqnarray}
K \hat{\mb{u}} & = & -M \hat{\mb{g}}.
\end{eqnarray}
The $\mb{\hat{u}}$ and $\mb{\hat{g}}$ denote the collection of coefficients $\hat{u}_j$ and $\hat{g}_j$ in the expansion of $\bar{u}$ and $\bar{g}$.  The $K$ denotes the stiffness matrix with entries $K_{ij} = \langle \mathcal{\tilde{L}} Y_j, Y_i \rangle_Q$ and $M$ denotes the mass matrix with entries $M_{ij} = \langle Y_j, Y_i \rangle_Q$.  We again mention that throughout the calculations we use expansions of functions up to order $\lfloor{L/2}\rfloor$ and Lebedev quadratures of order $L$. This is important to yield the needed over-sampling of functions to ensure accurate computation of the inner-product $\langle \cdot,\cdot \rangle_Q$.  

\section{Convergence Results}
\label{sec:Conv}
We discuss how our numerical methods converge in approximating the fundamental exterior calculus operations of the exterior derivative $\mb{d}$ and the Hodge star $\star$ when applied to $0$-forms, $1$-forms and $2$-forms.  We then discuss the convergence of our methods for compositions of operators and present results for the Laplace-Beltrami equation.  

We remark that throughout our convergence studies, we describe test functions using for a point $\mb{x}$ on the manifold its location within the embedding space.  To perform calculations we make use of the embedding space coordinates $[x,y,z]$ corresponding to $\mb{x} = x\bs{\iota}_1 + y\bs{\iota}_2  + z\bs{\iota}_3$, where $\bs{\iota}_1,\bs{\iota}_2,\bs{\iota}_3$ is the basis for the embedding space.  In this manner our test data is not tied to a specific choice of local coordinates on the manifold.  All figures report the relative error $\epsilon_{\mbox{\tiny rel}} = \|\mb{\bar{w}}^{\sharp} - \mb{w}^{\sharp}\|_2 / \|\mb{w}^{\sharp}\|_2$. The $\mb{w}$ is the exact result and $\mb{\bar{w}}$ is the numerically computed result. For $\|\cdot\|_2$, we use the $L^2$-norm of the embedding space.  

\subsection{Convergence of the Hodge Star Operator}
We approximate the Hodge star $\star$ by the numerical operator $\aHodgeStar$ using the hyperinterpolation approach discussed in Section~\ref{sec_num_ext}.  We investigate in practice the accuracy of this approach on a few different geometries and differential forms.

\begin{figure}[H]
\centering
\includegraphics[width=0.32\linewidth]{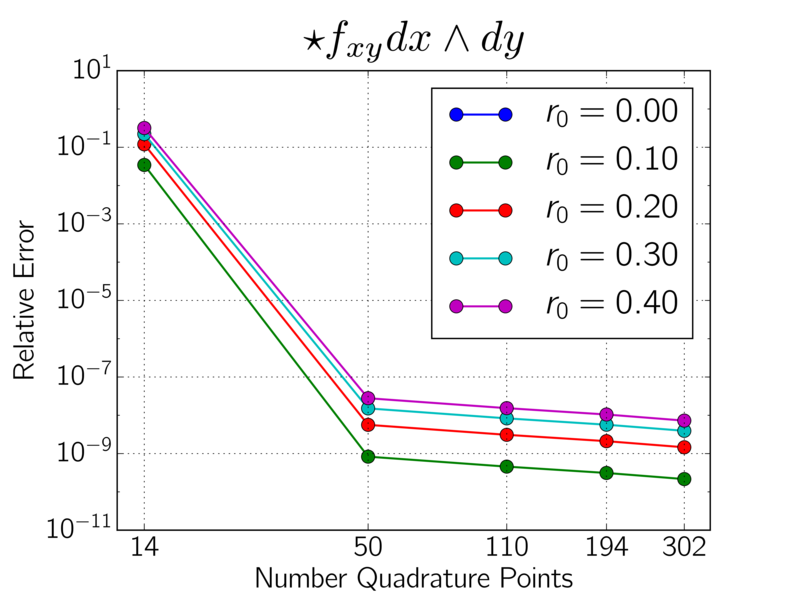}
\includegraphics[width=0.32\linewidth]{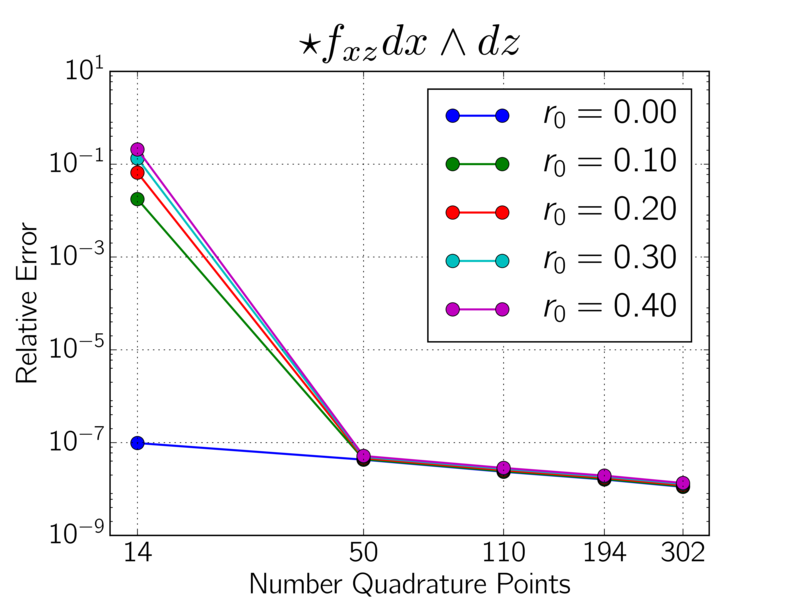}
\includegraphics[width=0.32\linewidth]{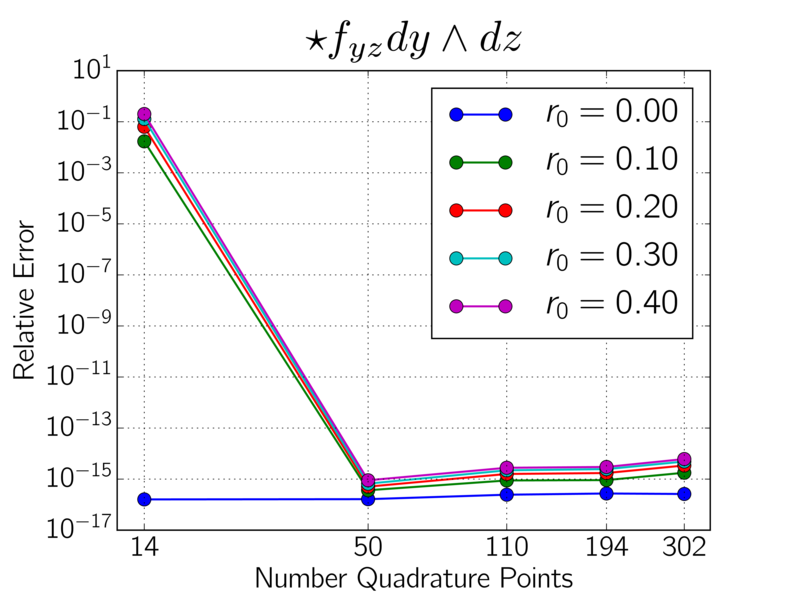}
\caption{Convergence of the numerical Hodge star operator $\aHodgeStar$ for $0$-forms.  We show for Manifold B how the relative error of  $\aHodgeStar f$ in approximating $\star f$ as the number of Lebedev nodes increases.  The $0$-form is $f = \exp(z)/(3 - y)$.  We investigate how the manifold geometry influences convergence by varying the amplitude $r_0$ in the range $[0.0,0.4]$ for Manifold B.  The amplitude $r_0 = 0.0$ corresponds to a sphere and $r_0 = 0.4$ to the final shape of Manifold B shown in Figure~\ref{fig:RadialManifolds}.}
\label{fig:HodgeStar_ZeroForm}
\end{figure}

We first consider a $0$-form defined on Manifold B defined in Figure~\ref{fig:RadialManifolds}.  We take the $0$-form to be $f = \exp(z)/(3 - y)$.  We show the accuracy of our numerical operator $\aHodgeStar$ in approximating the Hodge star $\star$ as the number of Lebedev nodes increases.  We find that the main limitation in the accuracy of the $\aHodgeStar$ is the resolution of the geometry of the manifold.  This is seen in our results where once a sufficient number of spherical harmonic modes are reached the relative error rapidly decays in approximating $\star f$.  The convergence of $\aHodgeStar$ as the number of Lebedev nodes is increased and when the geometry is varied is shown in Figure~\ref{fig:HodgeStar_ZeroForm}.

\begin{figure}[H]
\centering
\includegraphics[width=0.32\linewidth]{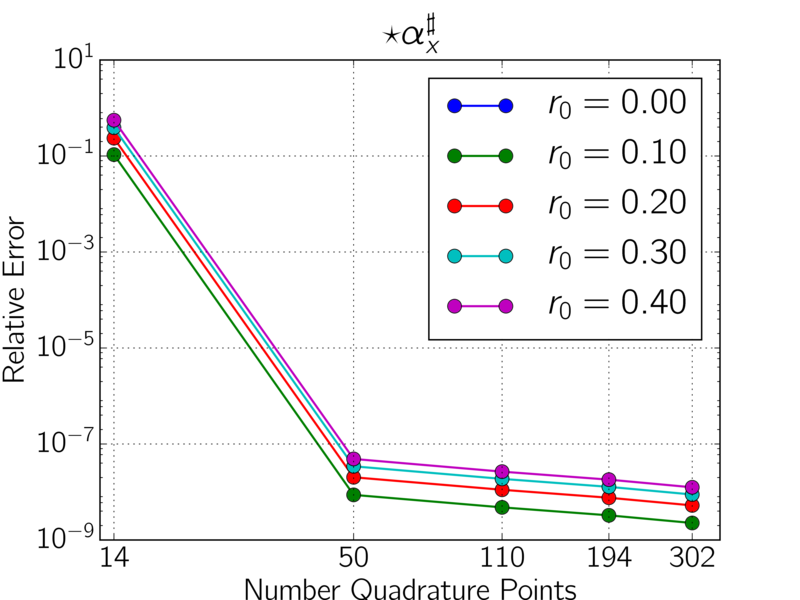}
\includegraphics[width=0.32\linewidth]{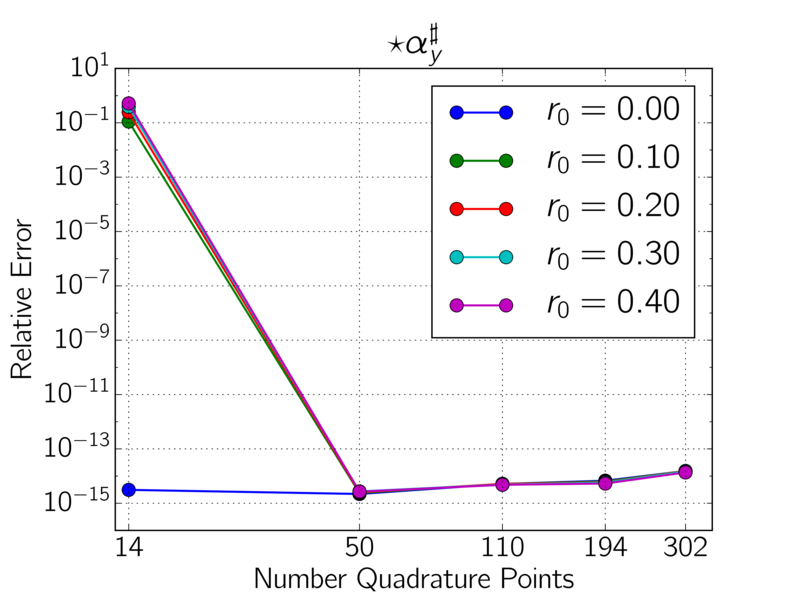}
\includegraphics[width=0.32\linewidth]{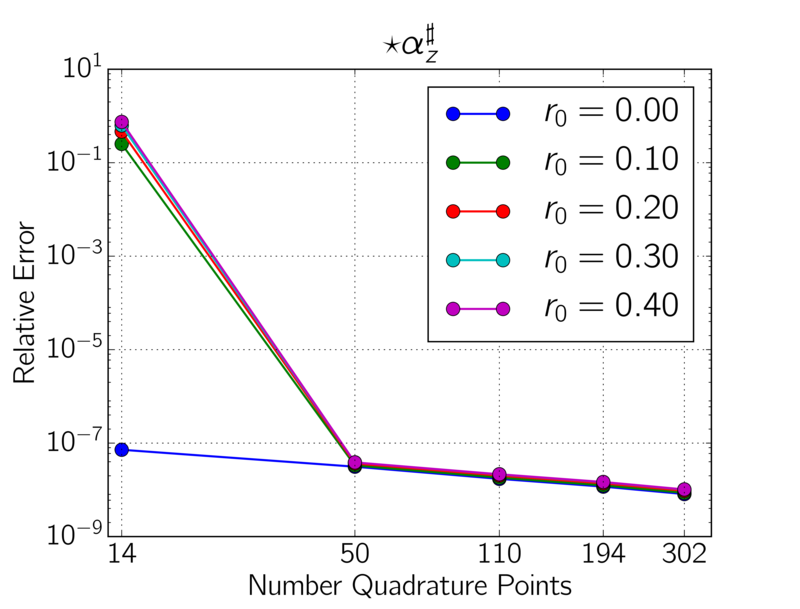}
\caption{Convergence of the numerical Hodge star operator $\aHodgeStar$ for $1$-forms.  We show for Manifold B the relative error of $\aHodgeStar\alpha$ in approximating $\star\alpha$ as the number of Lebedev nodes increases.  The $1$-form is $\alpha = \sqrt{|g|}\exp(z) d\theta + \sqrt{|g|}\exp(z) d\phi$.  We investigate how the manifold geometry influences convergence by varying the amplitude $r_0$ in the range $[0.0,0.4]$ for Manifold B.  The amplitude $r_0 = 0.0$ corresponds to a sphere and $r_0 = 0.4$ to the final shape of Manifold B shown in Figure~\ref{fig:RadialManifolds}.}
\label{fig:HodgeStar_OneForm}
\end{figure}

We next consider on Manifold B the $1$-form $\alpha = \sqrt{|g|}\exp(z) d\theta + \sqrt{|g|}\exp(z) d\phi$.  We again find that the main limitation in the accuracy of the $\aHodgeStar$ is the resolution of the geometry of the manifold.  In this case we find the error rapidly decreases once a sufficient number spherical harmonic modes are used.   The convergence of $\aHodgeStar$ as the number of Lebedev nodes is increased and when the geometry is varied is shown in Figure~\ref{fig:HodgeStar_OneForm}.

These results indicate that the numerical operator $\aHodgeStar$ provides an accurate approximation to the Hodge star $\star$.  Given the localized nature of the Hodge star, the main consideration to obtain accurate results with $\aHodgeStar$ is to use enough Lebedev nodes to ensure sufficient resolution of the geometry of the manifold.

\subsection{Convergence of the Exterior Derivative Operator}
\begin{figure}[H]
\centering
\includegraphics[width=0.32\linewidth]{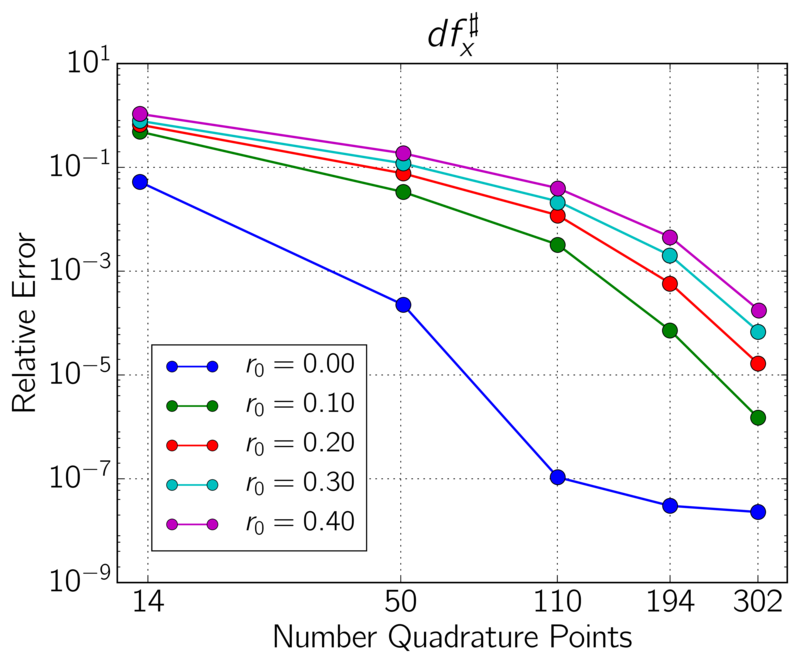}
\includegraphics[width=0.32\linewidth]{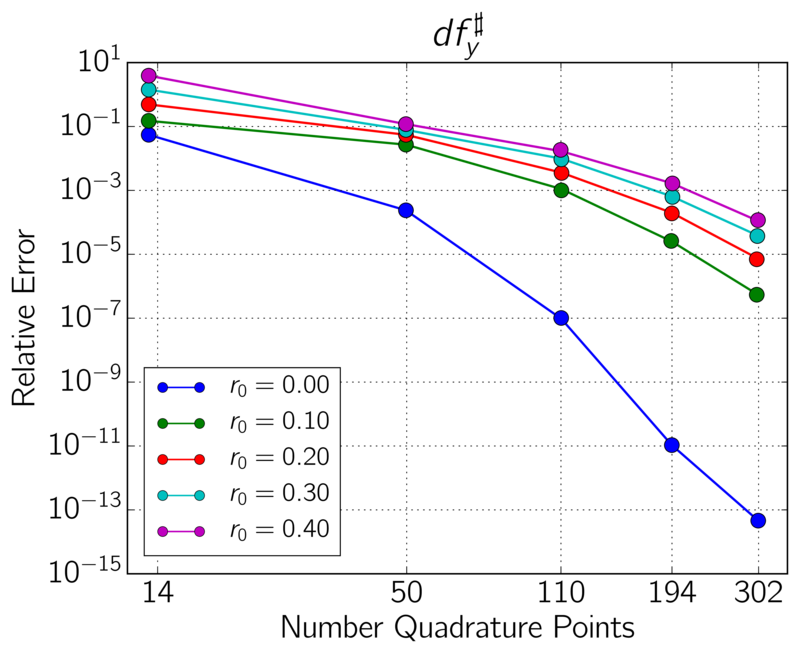}
\includegraphics[width=0.32\linewidth]{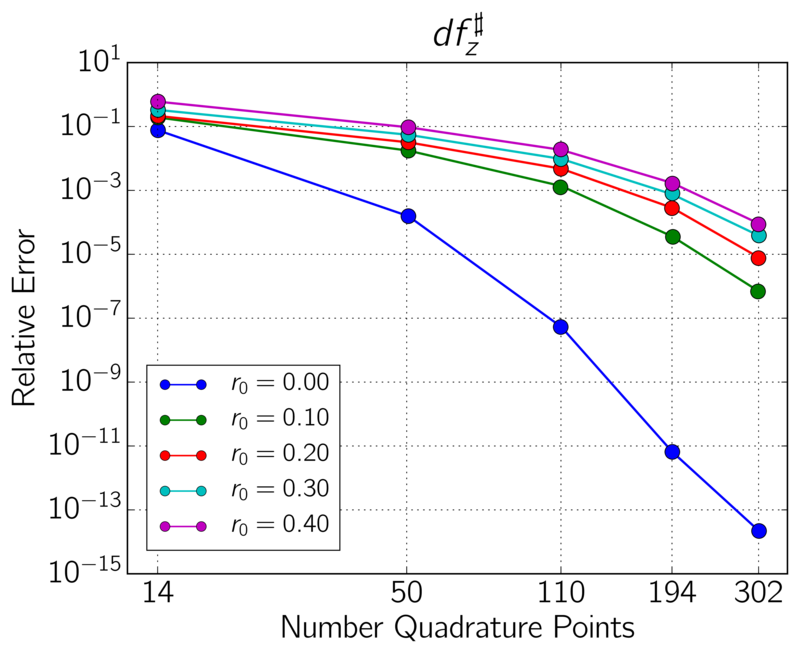}
\caption{Convergence of the numerical exterior derivative operator $\aExtD$ for $0$-forms.  We show for Manifold B the relative error of $\aExtD f$ in approximating $\mb{d}f$ as the number of Lebedev nodes increases.  The $0$-form is $f = \exp(z)$.  We investigate how the manifold geometry influences convergence by varying the amplitude $r_0$ in the range $[0.0,0.4]$ for Manifold B.  The amplitude $r_0 = 0.0$ corresponds to a sphere and $r_0 = 0.4$ to the final shape of Manifold B shown in Figure~\ref{fig:RadialManifolds}.}
\label{fig:ExteriorDerivativeTest1}
\end{figure}

\begin{figure}[H]
\centering
\includegraphics[width=0.32\linewidth]{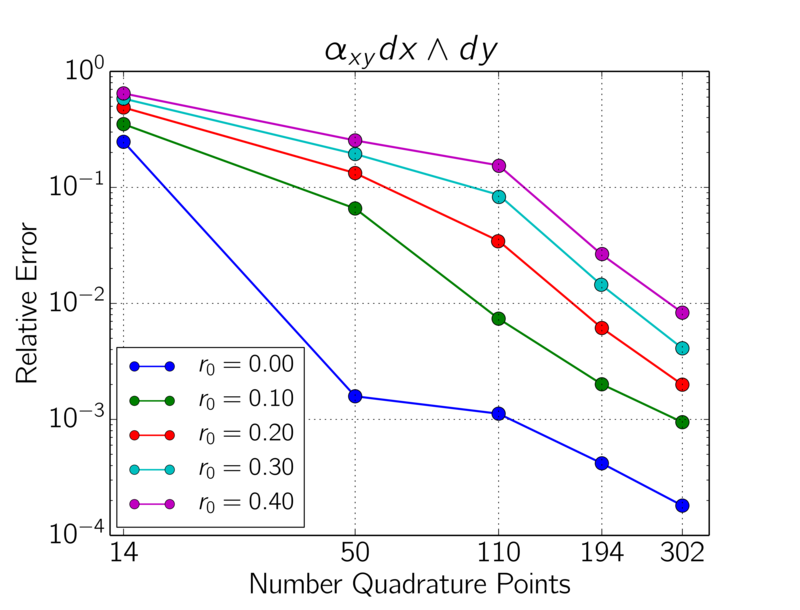}
\includegraphics[width=0.32\linewidth]{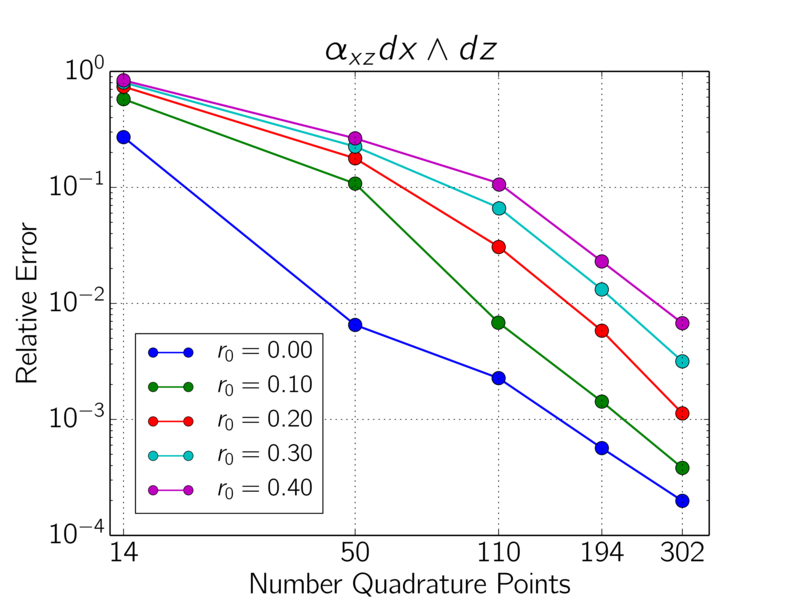}
\includegraphics[width=0.32\linewidth]{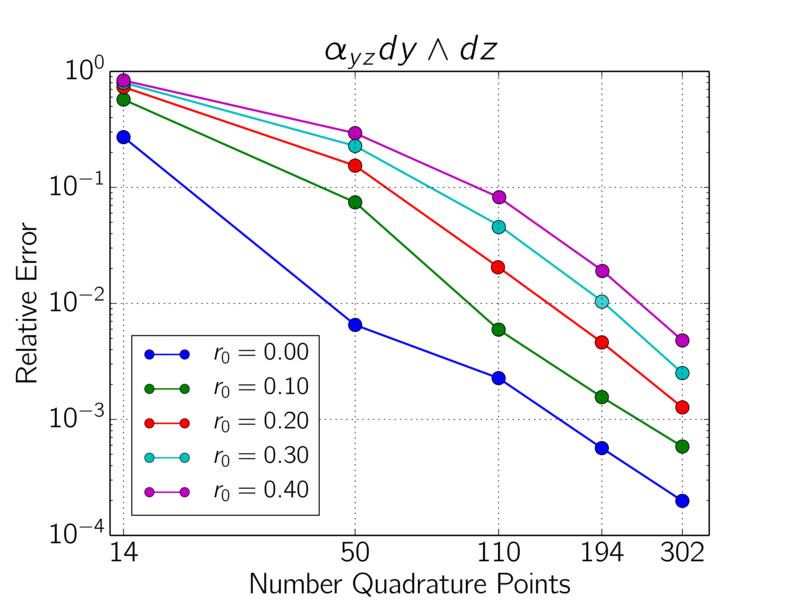}
\caption{Convergence of the numerical exterior derivative operator $\aExtD$ for $1$-forms.  We show for Manifold B the relative error of $\aExtD \alpha$ in approximating $\mb{d}\alpha$ as the number of Lebedev nodes increases.  The $1$-form is $\alpha = |g|\exp(z) d\theta + |g|\exp(z) d\phi$.  We investigate how the manifold geometry influences convergence by varying the amplitude $r_0$ in the range $[0.0,0.4]$ for Manifold B.  The amplitude $r_0 = 0.0$ corresponds to a sphere and $r_0 = 0.4$ to the final shape of Manifold B shown in Figure~\ref{fig:RadialManifolds}.}
\label{fig:ExteriorDerivativeTest2}
\end{figure}

We investigate the convergence of the numerical operator $\aExtD$ in approximating the exterior derivative $\mb{d}$ when applied to $0$-forms and $1$-forms.  We do not consider $2$-forms here since in two dimensions the exterior derivative would be zero which we also impose in our numerical calculations~\cite{Abraham1988,Spivak1971}.  We consider the $0$-form given by $f = \exp(z)$ and the $1$-form given by $\alpha = |g|\exp(z) d\theta + |g|\exp(z) d\phi$.  We consider for Manifold B the relative error of $\aExtD$ in approximating $\mb{d}$ as the number of Lebedev nodes is increased and the manifold geometry is varied.  These results are shown in Figure~\ref{fig:ExteriorDerivativeTest1} and Figure~\ref{fig:ExteriorDerivativeTest2}.

We see that the numerical operator $\aExtD$ converges spectrally in approximating the exterior derivative $\mb{d}$ both for the $0$-forms and $1$-forms.  Interestingly, we see that in the case when $r_0 = 0.0$ the approximation converges significantly more rapid than the cases when $r_0 \neq 0$.  This occurs since $r_0 = 0$ corresponds to the case when the shape is a sphere where the geometry is relatively simple and many of the geometric terms to be numerically approximated greatly simplify.  The main source of error in this case arises primarily from the hyperinterpolation used for computing the derivatives.  In the case when $r_0 \neq 0$, the isogeometric approach used to compute $\aExtD$ approximates the geometry of the manifold using a finite spherical harmonics representation which results in an additional source of approximation error.  

We investigate how the geometry contributes to convergence by varying $r_0$ over the range $[0.0,0.4]$.  The case with $r_0 = 0.4$ corresponds to the final shape of Manifold B shown in Figure~\ref{fig:RadialManifolds}.  The convergence is found to be spectral and comparable in each of these different cases for the geometry.  We see that as one might expect the largest errors are incurred in the case with the most pronounced geometry corresponding to $r_0 = 0.4$.  Overall, the results indicate that the numerical operator $\aExtD$ provides for $0$-forms and $1$-forms an accurate approximation for the exterior derivative $\mb{d}$.

\subsection{Convergence for the Laplace-Beltrami Equation}

We consider the composition of exterior calculus operators to represent partial differential equations on manifolds.  Taking the approach we discussed in Section~\ref{sec:methodsForPDE}, we consider the Laplace-Beltrami operator $\mathcal{L} = -\bs{\delta} \mb{d} = -\star \mb{d} \star \mb{d}$ and consider the Poisson problem on the manifold 
\begin{eqnarray}
\label{equ:PoissonEqu}
\mathcal{L}u = -\bs{\delta} \mb{d} u = -g.
\end{eqnarray}
Since there are no boundaries for the manifolds we shall consider, we also impose throughout that the zero mode has $\hat{u}_0 = 0$.  The Laplace-Beltrami equation~\ref{equ:PoissonEqu} provides a test of using compositions of the numerical exterior calculus operators $\aExtD$ and $\aHodgeStar$ and the associated convergence and accuracy.  

\begin{figure}[H]
\centering
\includegraphics[width=0.8\linewidth]{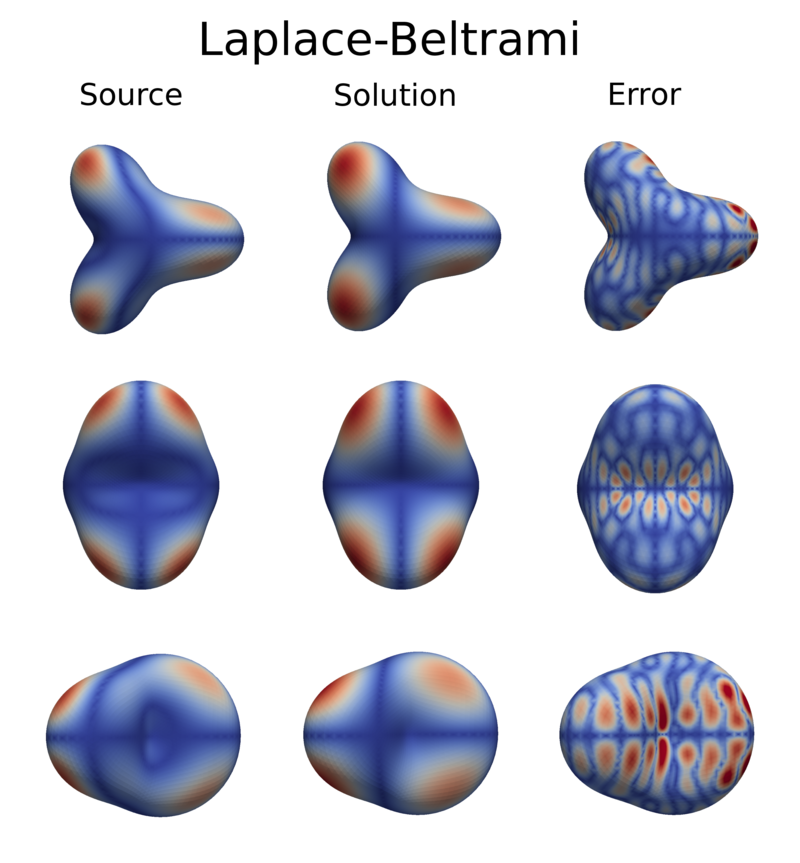}
\caption{Solution of the Laplace-Beltrami Equations on Manifold B of Figure~\ref{fig:RadialManifolds}.  For the Laplace-Beltrami equations we show on the manifold surface the source term $g$ (left), solution $u$ (middle), and relative error $\epsilon_{rel}$ (right) from the vantage point of the coordinate axes associated with the embedding space $\bs{\iota}_1$, $\bs{\iota}_2$, $\bs{\iota}_3$.  Shown is the case with the spherical harmonics resolved with $432$ Lebedev nodes.}    
\label{fig:DimpleSolutionViewsEach}
\end{figure}

We investigate how our methods solve the Laplace-Beltrami equation~\ref{equ:PoissonEqu} on Manifold B and Manifold C of Figure~\ref{fig:RadialManifolds}.  To have a known solution to the Laplace-Beltrami equations with which to compare, we manifacture a source term for the solution on the manifold given by $u = \exp(y)/(3-z)^4$, where the $x,y,z$ refer to the coordinates of the embedding space as discussed in the beginning of Section~\ref{sec:Conv}.  We compute the source term $g= -\mathcal{L}u$ symbolically using the package SymPy and evaluate the expressions numerically when data is needed for the source~\cite{Sympy2017}.  In this manner we are able to assess the relative errors of the numerical methods with high-precision.

\begin{figure}[H]
\centering
\includegraphics[width=0.9\linewidth]{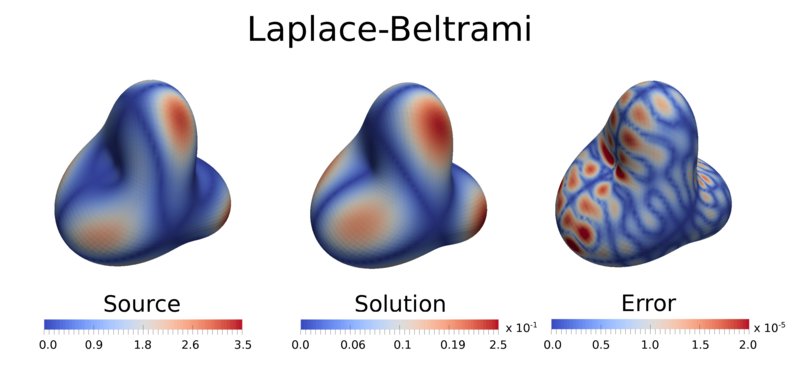}
\caption{Solution of the Laplace-Beltrami Equations on Manifold B of Figure~\ref{fig:RadialManifolds}.  For the Laplace-Beltrami equations we show on the manifold surface the source term $g$ (left), solution $u$ (middle), and relative error $\epsilon_{rel}$ (right) from a rotated vantage point.  The case shown corresponds to $r_0 = 0.4$ and spherical harmonics resolved with $432$ Lebedev nodes.}
\label{fig:DimpleSolutionViewDiag}
\end{figure}

\begin{figure}[H]
\centering
\includegraphics[width=0.6\linewidth]{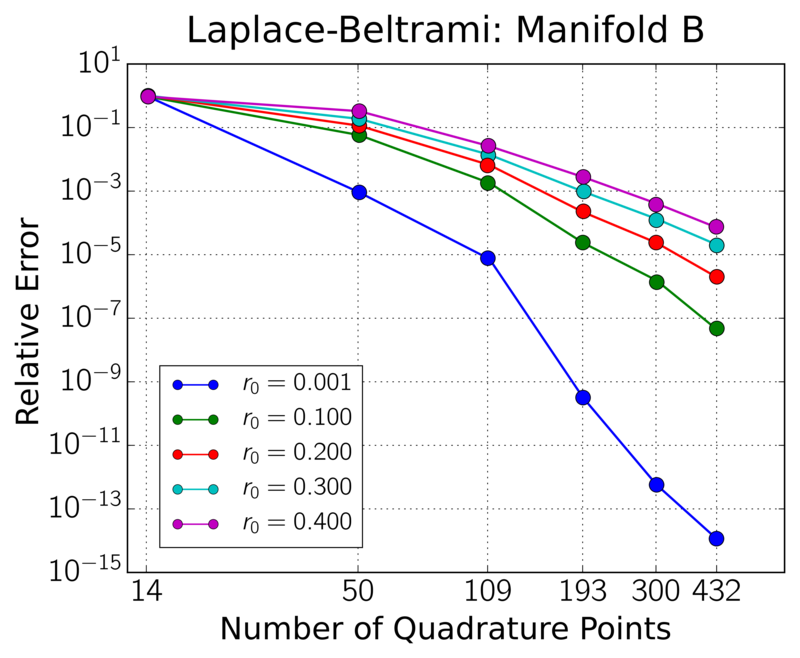}
\caption{Convergence of the numerical methods for the Laplace-Beltrami Equations on Manifold B of Figure~\ref{fig:RadialManifolds}.  We show for Manifold B the relative error when the Laplace-Beltrami equations are solved using a composition of $\aExtD$ and $\aHodgeStar$ as discussed in Section~\ref{sec:methodsForPDE}.  We investigate how the manifold geometry influences convergence by varying the amplitude $r_0$ in the range $[0.0,0.4]$ for Manifold B.  The amplitude $r_0 = 0.0$ corresponds to a sphere and $r_0 = 0.4$ to the final shape of Manifold B shown in Figure~\ref{fig:RadialManifolds}.}
\label{fig:Dimple_LB_Conv}
\end{figure}

We show the source term $g$, solution $u$, and relative errors $\epsilon_{rel}$ in Figure~\ref{fig:DimpleSolutionViewsEach} and Figure~\ref{fig:DimpleSolutionViewDiag}.  We use the approach discussed for solving PDEs on manifolds in Section~\ref{sec:methodsForPDE} for the Laplace-Beltrami equations.  We show how the composition of the numerical operators $\aExtD$ and $\aHodgeStar$ perform as the number of Lebedev nodes increases in Figure~\ref{fig:Dimple_LB_Conv}. The results indicate that the numerical methods based on composing $\aExtD$ and $\aHodgeStar$ perform well in approximating the true compositions of the exterior calculus operators $\star$ and $\mb{d}$.

We see that when the geometry has $r_0 = 0.0$ corresponding to a sphere the methods converge most rapidly.  In the cases with $r_0 \neq 0$ we see the convergence is comparable in each case. We see the convergence is slowest for the most pronounced geometry where the representation for the geometry incurs the most approximation error.  Overall, we find the compositions of $\aExtD$ and $\aHodgeStar$ perform well and provide a solver with spectral accuracy for the Laplace-Beltrami equations on Manifold B.     

\begin{figure}[H]
\centering
\includegraphics[width=0.8\linewidth]{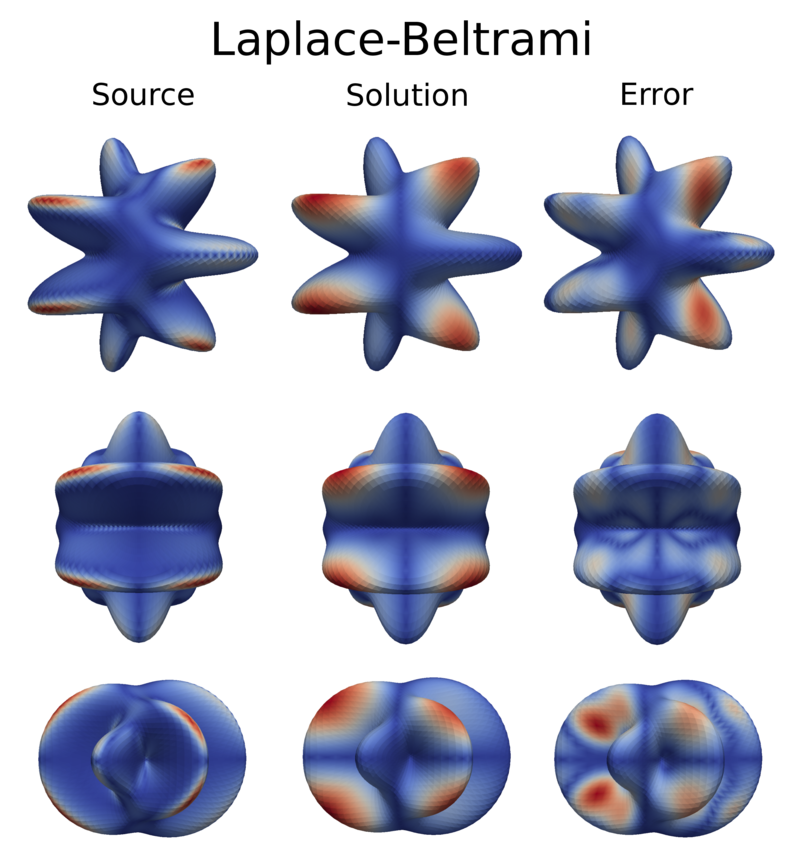}
\caption{Solution of the Laplace-Beltrami Equations on Manifold C of Figure~\ref{fig:RadialManifolds}.  For the Laplace-Beltrami equations we show on the manifold surface the source term $g$ (left), solution $u$ (middle), and relative error $\epsilon_{rel}$ (right) from the vantage point of the coordinate axes associated with the embedding space $\bs{\iota}_1$, $\bs{\iota}_2$, $\bs{\iota}_3$.  Shown is the case with the spherical harmonics resolved with $302$ Lebedev nodes.}
\label{fig:FountainSolutionViewsEach}
\end{figure}

\begin{figure}[H]
\centering
\includegraphics[width=0.9\linewidth]{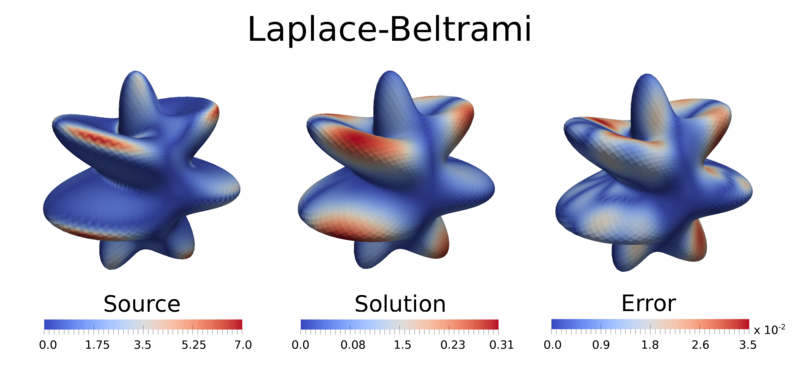}
\caption{Solution of the Laplace-Beltrami Equations on Manifold C of Figure~\ref{fig:RadialManifolds}.  For the Laplace-Beltrami equations we show on the manifold surface the source term $g$ (left), solution $u$ (middle), and relative error $\epsilon_{rel}$ (right) from a rotated vantage point.  The case shown corresponds to $r_0 = 0.4$ and the spherical harmonics resolved with $302$ Lebedev nodes.}
\label{fig:FountainSolutionViewsDiag}
\end{figure}

We next consider the more geometrically complex Manifold C, see Figure~\ref{fig:RadialManifolds}.  We again use the manufactured solution on the surface given by $u = \exp(y)/(3-z)^4$ and $g = -\mathcal{L}u$ computed symbolically.  We mention that given the different geometry of Manifold C relative to Manifold B both the solution and the source data are different on each of the manifold surfaces.  We show on the manifold surface the source term $g$, solution $u$, and relative error $\epsilon_{rel}$ in Figure~\ref{fig:FountainSolutionViewsEach} and Figure~\ref{fig:FountainSolutionViewsDiag}.  

\begin{figure}[H]
\centering
\includegraphics[width=0.6\linewidth]{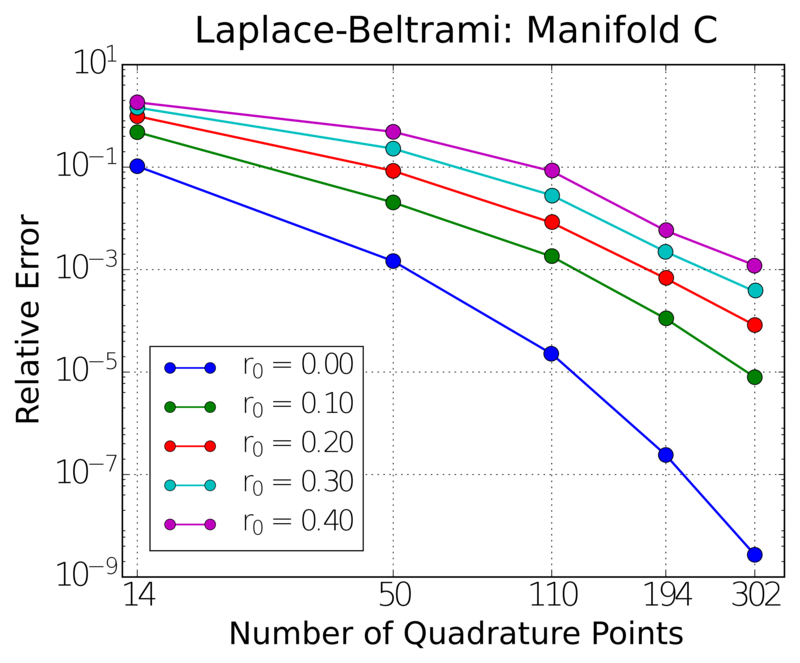}
\caption{Convergence of the numerical methods for the Laplace-Beltrami Equations on Manifold C of Figure~\ref{fig:RadialManifolds}.  We show for Manifold B the relative error when the Laplace-Beltrami equations are solved using a composition of $\aExtD$ and $\aHodgeStar$ as discussed in Section~\ref{sec:methodsForPDE}.  We investigate how the manifold geometry influences convergence by varying the amplitude $r_0$ in the range $[0.0,0.4]$ for Manifold C.  The amplitude $r_0 = 0.0$ corresponds to a sphere and $r_0 = 0.4$ to the final shape of Manifold C shown in Figure~\ref{fig:RadialManifolds}.}
\label{fig:Fountain_LB_Conv}
\end{figure}

For the more geometrically complicated Manifold C, we perform calculations to show how the composition of the numerical operators $\aExtD$ and $\aHodgeStar$ perform as the number of Lebedev nodes increases and $r_0$ is varied, see Figure~\ref{fig:Fountain_LB_Conv}. 

We find the numerical methods converge spectrally for Manifold C.  As may be expected given its more complicated geometry, we find the convergence is somewhat slower than for Manifold B.  We find that the numerical methods perform well for Manifold C once a sufficient number of Lebedev nodes are used to resolve the geometry.  We again see the trend that the case with smaller $r_0$ convergence more rapidly.  Overall, we find that in practice the composition of the numerical operator $\aExtD$ and $\aHodgeStar$ performs well in approximating differential operators on the manifold.  This indicates the approach discussed in Section~\ref{sec:methodsForPDE} can be applied quite generally in practice for solving partial differential equations on manifolds.  

\section{Conclusion}
In summary, we have developed numerical methods for exterior calculus on radial manifolds based on hyperinterpolation and Lebedev quadrature.    We have shown our methods are spectrally accurate in approximating the exterior derivative $\mb{d}$, Hodge star $\star$, and their compositions.  We presented results for the Laplace-Beltrami equations that demonstrate the utility of the approach.  The introduced numerical methods can be applied quite generally for approximating exterior calculus operations and differential equations on radial manifolds.

\section{Acknowledgments}
We acknowledge support to P.J.A. and B.G. from research grants DOE ASCR CM4 DE-SC0009254, NSF CAREER Grant DMS-0956210, and NSF Grant DMS - 1616353.

\bibliographystyle{plain}
\bibliography{paperDatabase}{}

\appendix

\section*{Appendix}

\section{Spherical Harmonics}
\label{appendix:sphericalHarmonics}
The spherical harmonics are given by
\begin{equation}
\label{equ:sphericalHarmonics}
Y^m_n(\theta, \phi) = \sqrt{\frac{(2n+1)(n-m)!}{4\pi(n+m)!}}P^m_n\left(\cos(\phi)\right) \exp\left({im\theta}\right)
\end{equation}
where $m$ denotes the order and $n$ the degree for $n \ge 0$ and $m \in \{-n, \dots, n\}$.  The $P^m_n$ denote the \textit{Associated Legendre Polynomials}.  In our notation, $\theta$ denotes the azmuthal angle and $\phi$ the polar angle of the spherical coordinates~\cite{HanBookSphericalHarmonics2010}.

Since we work throughout only with real-valued functions, we have that the modes are self-conjugate and we use that $Y^m_n = \overline{Y^{-m}_n}$.  We have found it convenient to represent the spherical harmonics as
\begin{equation}
Y^m_n(\theta, \phi) = X^m_n(\theta, \phi)+ iZ^m_n(\theta, \phi) 
\end{equation}
where $X_n^m$ and $Z_n^m$ denote the real and imaginary parts.  In our numerical methods we use this splitting to construct a purely real set of basis functions on the unit sphere with maximum degree $N$.  We remark that this consists of $(N+1)^2$ basis elements.  In the case of $N=2$ we have the basis elements  
\begin{eqnarray}
\tilde{Y}_1 = Y^0_0,\hspace{0.2cm} 
\tilde{Y}_2 = Z^1_1,\hspace{0.2cm} 
\tilde{Y}_3 = Y^0_1,\hspace{0.2cm} 
\tilde{Y}_4 = X^1_1,\hspace{0.2cm} 
\tilde{Y}_5 = Z^2_2,\hspace{0.2cm} 
\tilde{Y}_6 = Z^1_2,\hspace{0.2cm} 
\tilde{Y}_7 = Y^0_2,\hspace{0.2cm} 
\tilde{Y}_8 = X^1_2,\hspace{0.2cm} 
\tilde{Y}_9 = X^2_2.\hspace{0.2cm} 
\end{eqnarray}
We use a similar convention for the basis for the other values of $N$.  We take our final basis elements $Y_i$ to be the normalized as $Y_i = \tilde{Y}_i/\sqrt{\langle \tilde{Y}_i, \tilde{Y}_i \rangle}$.

We compute derivatives of our finite expansions by evaluating analytic formulas for the spherical harmonics in order to try to minimize approximation error~\cite{HanBookSphericalHarmonics2010}.  Approximation errors are incurred when sampling the values of expressions involving these derivatives at the Lebedev nodes and performing quadratures.  The derivative in the azmuthal coordinate $\theta$ of the spherical harmonics is given by
\begin{equation*}
\partial_\theta Y^m_n(\theta, \phi) 
= \partial_\theta \sqrt{\frac{(2n+1)(n-m)!}{4\pi(n+m)!}}P^m_n(\cos(\phi)) \exp\left({im\theta}\right)
= imY^m_n\left(\theta, \phi\right).
\end{equation*}
This maps the spherical harmonic of degree $n$ to again a spherical harmonic of degree $n$.  In our numerics, this derivative can be represented in our finite basis which allows us to avoid projections.  This allows for computing the derivative in $\theta$ without incurring an approximation error.  For the derivative in the polar angle $\phi$ we have that
\begin{equation}
\label{equ:derivPhi}
\partial_\phi Y^m_n(\theta, \phi) 
= m \cot(\phi)Y^m_n(\theta, \phi)+\sqrt{(n-m)(n+m+1)}\exp\left({-i\theta}\right) Y^{m+1}_n(\theta, \phi).
\end{equation}
We remark that the expression can not be represented in terms of a finite expansion of spherical harmonics.  We use this expression for $\partial_\phi Y^m_n(\theta, \phi)$ when we compute 
values at the Lebedev quadrature nodes in equation~\ref{equ_der_comp_v}.  This provides a convenient way to compute derivatives of differential forms following the approach discussed in Section~\ref{sec:num_methods}.  We remark that it is the subsequent hyperinterpolation of the resulting expressions where the approximation error is incurred.  We adopt the notational convention that $Y^{m}_n = 0$ when $m \geq n+1$.  For further discussion of spherical harmonics see~\cite{HanBookSphericalHarmonics2010}.

\section{Differential Geometry of Radial Manifolds}
\label{appendix:radialManifold}
We consider throughout manifolds of radial shape.  A radial manifold is defined as a surface where each point can be connected by a line segment to the origin without intersecting the surface.  In spherical coordinates, any point $\mb{x}$ on the radial manifold can be expressed as
\begin{eqnarray}
\label{equ:manifoldParam}
\mb{x}(\theta,\phi) = \bs{\sigma}(\theta, \phi) =  r(\theta, \phi)\mb{r}(\theta,\phi)
\end{eqnarray}
where $\mb{r}$ is the unit vector from the origin to the point on the sphere corresponding to angle $\theta,\phi$ and $r$ is a positive scalar function.

We take an isogeometric approach to representing the manifold $M$.  We sample the scalar function $r$ at the Lebedev nodes and represent the geometry using the finite spherical harmonics expansion $r(\theta,\phi) = \sum_i \bar{r}_i Y_i$ up to the order $\lfloor L/2 \rfloor$ where $\bar{r}_i = \langle r, Y_i \rangle_Q$ for a quadrature of order $L$.

We consider two coordinate charts for our calculations.  The first is referred to as Chart A and has coordinate singularities at the north and south pole.  The second is referred to as Chart B and has coordinate singularities at the east and west pole.  For each chart we use spherical coordinates with $(\theta, \phi) \in [0 , 2\pi) \times [0 , \pi]$ but to avoid singularities only use values in the restricted range $\phi \in [\phi_{min}, \phi_{max}]$, where $ 0 < \phi_{min} \le \frac{\pi}{4}$, and $\frac{3\pi}{4} \le \phi_{max} < \pi$.  In practice, one typically takes $\phi_{min} = 0.8\times \frac{\pi}{4}$ and $\phi_{max} = 0.8 \times \pi$. 
For Chart A, the manifold is parameterized in the embedding space $\mathbb{R}^3$ as
\begin{eqnarray}
\mb{x}(\hat{\theta}, \hat{\phi}) = r(\hat{\theta}, \hat{\phi})\mathbf{r}(\hat{\theta}, \hat{\phi}), \hspace{0.5cm} 
\mathbf{r}(\hat{\theta}, \hat{\phi}) = \begin{bmatrix}\sin(\hat{\phi})\cos(\hat{\theta}), & \sin(\hat{\phi})\sin(\hat{\theta}), & \cos(\hat{\phi}) \end{bmatrix} 
\label{equ:defChartA}
\end{eqnarray}
and for Chart B 
\begin{eqnarray}
\mb{x}(\bar{\theta}, \bar{\phi}) = r(\bar{\theta}, \bar{\phi})\mathbf{r}(\bar{\theta}, \bar{\phi}), \hspace{0.5cm} 
\mathbf{\bar{r}}(\bar{\theta}, \bar{\phi}) = \begin{bmatrix} \cos(\bar{\phi}), & \sin(\bar{\phi})\sin(\bar{\theta}), & -\sin(\bar{\phi})\cos(\bar{\theta}) \end{bmatrix}.
\label{equ:defChartB}
\end{eqnarray}
With these coordinate representations, we can derive explicit expressions for geometric quantities associated with the manifold such as the metric tensor and shape tensor.  The derivatives used as the basis $\partial_\theta, \partial_\phi$ for the tangent space can be expressed as
\begin{eqnarray}
\label{eqn:sigmaTheta}
\bs{\sigma}_{\theta}(\theta, \phi) & = & r_{\theta}(\theta, \phi)\mb{r}(\theta, \phi) + r(\theta, \phi) \mb{r}_{\theta}(\theta, \phi)\\
\label{eqn:sigmaPhi}
\bs{\sigma}_{\phi}(\theta, \phi) & = & r_{\phi}(\theta, \phi)\mb{r}(\theta, \phi) + r(\theta, \phi) \mb{r}_{\phi}(\theta, \phi).
\end{eqnarray}
We have expressions for $\mb{r}_{\theta}$ and $\mb{r}_{\phi}$ in the embedding space $\mathbb{R}^3$ using equation~\ref{equ:defChartA} or equation~\ref{equ:defChartB} depending on the chart being used.  The first fundamental form $\mathbf{I}$ (metric tensor) and second fundamental form $\mathbf{II}$ (shape tensor) are given by 
\begin{align}
\mathbf{I} = \begin{bmatrix}
E & F \\
F & G
\end{bmatrix} = \begin{bmatrix}
\bs{\sigma}_{\theta} \cdot \bs{\sigma}_{\theta} & \bs{\sigma}_{\theta} \cdot \bs{\sigma}_{\phi} \\
\bs{\sigma}_{\phi} \cdot \bs{\sigma}_{\theta} & \bs{\sigma}_{\phi} \cdot \bs{\sigma}_{\phi}
\end{bmatrix} =  \begin{bmatrix}
r_{\theta}^2+r^2\sin(\phi)^2 & r_{\theta}r_{\phi} \\
r_{\theta}r_{\phi} & r_{\phi}^2+r^2
\end{bmatrix}. \label{equ:I_mat}
\end{align}
and 
\begin{align}
\mathbf{II} = \begin{bmatrix}
L & M \\
N & N
\end{bmatrix} = \begin{bmatrix}
\bs{\sigma}_{\theta \theta} \cdot \bs{n} & \bs{\sigma}_{\theta \phi} \cdot \bs{n} \\
\bs{\sigma}_{\phi \theta} \cdot \bs{n} & \bs{\sigma}_{\phi \phi}\cdot \bs{n}
\end{bmatrix}. \label{equ:II_mat} 
\end{align}
The $\mb{n}$ denotes the outward normal on the surface and is computed using
\begin{eqnarray}
\bs{n}(\theta, \phi) = \frac{\bs{\sigma}_{\theta}(\theta, \phi) \times \bs{\sigma}_{\phi}(\theta, \phi)}{\| \bs{\sigma}_{\theta}(\theta, \phi) \times \bs{\sigma}_{\phi}(\theta, \phi) \|}. \label{equ:normal_vec_def}
\end{eqnarray}
The terms $\bs{\sigma}_{\theta\theta}$, $\bs{\sigma}_{\theta\phi}$, and $\bs{\sigma}_{\phi,\phi}$ are obtained by further differentiation from equation~\ref{eqn:sigmaTheta} and equation~\ref{eqn:sigmaPhi}.
We use the notation for the metric tensor $\mb{g} = \mathbf{I}$  interchangably.  In practical calculations whenever we need to compute the action of the inverse metric tensor we do so through numerical linear algebra (Gaussian elimination with pivoting)~\cite{Trefethen1997,Strang1980}.  For notational convenience, we use the tensor notation for the metric tensor $g_{ij}$ and its inverse $g^{ij}$ which has the formal correspondence
\begin{eqnarray}
g_{ij} = \left[\mathbf{I}\right]_{i,j}, \hspace{0.3cm} g^{ij} = \left[ \mathbf{I}^{-1}\right]_{i,j}.
\end{eqnarray}
For the metric factor we also have that
\begin{eqnarray}
\sqrt{|g|} = \sqrt{\det(\mathbf{I})} = r\sqrt{r_{\theta}^2+(r_{\phi}^2+r^2)\sin(\phi)^2} = \| \vec{\sigma}_{\theta}(\theta, \phi) \times \vec{\sigma}_{\phi}(\theta, \phi) \|.
\label{met_fac_def}
\end{eqnarray}
To ensure accurate numerical calculations in each of the above expressions the appropriate coordinates either Chart A or Chart B are used to ensure sufficient distance from coordinate singularities at the poles.  To compute quantities associated with curvature of the manifold we use the Weingarten map~\cite{Pressley2001} which can be expressed as
\begin{eqnarray}
\mb{W} = -\mb{I}^{-1} \mb{II}.
\end{eqnarray}
To compute the Gaussian curvature $K$, we use
\begin{eqnarray}
K(\theta,\phi) = \det\left(\mb{W}(\theta,\phi)\right).
\end{eqnarray}
For further discussion of the differential geometry of manifolds see~\cite{Pressley2001,Abraham1988,Spivak1971}.

\end{document}